\documentclass{article}
\usepackage{mathtools}
\usepackage{amssymb}
\usepackage[margin=2.325cm]{geometry}
\usepackage{titling}

\newcommand{\ax}{\mathcal{A}_x}
\newcommand{\hax}{\hat{\mathcal{A}}_x}
\newcommand{\abx}{\mathcal{A}_{\bar{x}}}
\newcommand{\uhax}{U(\hax)}
\newcommand{\uabx}{U(\abx)}
\newcommand{\lax}{\ell(\ax)}
\newcommand{\lhax}{\ell(\hax)}
\newcommand{\lanx}{\ell(\abx)}

\begin{document}

\title{\huge An upper bound for union-closed family size}
\author{\Large Christopher Bouchard}
\date{}
\maketitle

\vspace{-0.675cm}

\abstract{\vspace{-0.0675cm} \noindent Let $\mathcal{A}$ be a union-closed family of sets with universe $\bigcup_{A \in \mathcal{A}}A = [n] = \{1,\cdots,n\}$ and length $\ell$. We prove that $|\mathcal{A}| \leq \sum_{i=0}^{\ell} \binom{n}{i}$, with equality if and only if $\mathcal{A} = \bigcup_{i=0}^{\ell}\binom{[n]}{n-i}$. Additionally, by showing that $|\mathcal{A}| \leq \frac{\ell^p-1}{\ell-1}+2^n(1-2^{-\ell})^p$ for any nonnegative integer $p$, we establish for all integers $1 \leq k \leq n$ that $\sum_{i=0}^k \binom{n}{i} \leq \frac{k^{\hat{p}}-1}{k-1}+2^n(1-2^{-k})^{\hat{p}}$, where $\hat{p}=\lfloor (n-k)/\log_2(\frac{k}{1-2^{-k}})\rfloor + 1$.}

\vspace{0.675cm}

\section*{1. Introduction}

Let $\mathcal{A}$ be a finite family of distinct finite sets (at least one of which is nonempty) with \textit{universe} $U(\mathcal{A}) \coloneqq \bigcup_{A \in \mathcal{A}}A$ denoted by $[n]=\{1,\cdots,n\}$. $\mathcal{A}$ is called \textit{union-closed} if $X_1,X_2 \in \mathcal{A}$ implies that $X_1 \cup X_2 \in \mathcal{A}$. Such families have gained popularity by way of the union-closed sets conjecture, also known as Frankl's conjecture, which states that in any union-closed $\mathcal{A}$, there must be an element from $[n]$ that appears in at least half of its member sets (see [1--2] for an overview of many related results). One result, proved by Reimer in [4], can be expressed as an implicit least upper bound for the size of $\mathcal{A}$, namely $|\mathcal{A}| \leq (2 / \log_2|\mathcal{A}|)\sum_{A \in \mathcal{A}}|A|$ (or $|\mathcal{A}| \leq 4^{\sum_{A \in \mathcal{A}}|A|/|\mathcal{A}|}$).

A \textit{chain} $\mathcal{C}$ in $\mathcal{A}$ is a subfamily of $\mathcal{A}$ such that $X_1 , X_2 \in \mathcal{C}$ implies that $(X_1 \subseteq X_2) \lor (X_2 \subseteq X_1)$, and the \textit{length} of $\mathcal{A}$, denoted by $\ell \coloneqq \ell(\mathcal{A})$, is one less than the maximum size of a chain in $\mathcal{A}$. A result of Erd\H os (see Theorem 5 of [3]) states that any family of sets $\mathcal{A}$ (not necessarily union-closed) has size less than or equal to the sum of the largest $\ell+1$ binomial coefficients of $n$. In the present work, we prove that imposing the union-closed constraint on $\mathcal{A}$ tightens this upper bound to be the sum of the first, rather than largest, $\ell+1$ binomial coefficients. In a similar vein, we also provide an upper bound for the sum of the first $k$ binomial coefficients of any positive integer $n$. Denote by $\binom{S}{r}$ the family of $r$-subsets of a set $S$. The two main theorems are stated as follows:

\bigskip
\smallskip

\noindent \textbf{Theorem 1.} \textit{For any union-closed family} $\mathcal{A}$ \textit{with universe} $[n]$ \textit{and length} $\ell$\textit{,}\[|\mathcal{A}| \leq \sum_{i=0}^{\ell} \binom{n}{i}\textrm{,}\]

\vspace{0.375cm}

\noindent \textit{with equality and if and only if}
\vspace{-0.65cm}
\[\hspace{0.115cm} \mathcal{A} = \bigcup_{i=0}^{\ell}\binom{[n]}{n-i}\textrm{.}\]

\bigskip
\medskip

\noindent \textbf{Theorem 2.} \textit{For all integers} $1 \leq k \leq n$\textit{,}
\smallskip
\[\sum_{i=0}^k\binom{n}{i} \leq \frac{k^{\hat{p}}-1}{k-1}+2^n(1-2^{-k})^{\hat{p}}\textrm{,}\]

\vspace{0.375cm}

\noindent \textit{where}
\vspace{-0.5cm}
\[\hspace{0.225cm} \hat{p} = \Bigr \lfloor \frac{n-k}{\log_2(k/(1-2^{-k}))} \Bigr\rfloor + 1\textrm{.}\]

\medskip

\section*{2. Proof of Theorem 1}

We prove the theorem by induction on universe size. For the base case $n=1$, the only union-closed families are $\mathcal{A}=\{1\}$ and $\mathcal{A}=\{\{1\},\emptyset\}$. The former has $\ell=0$ and $|\mathcal{A}|=\sum_{i=0}^{\ell} \binom{n}{i} = \binom{1}{0}=1$, while the latter has $\ell=1$ and $|\mathcal{A}| = \sum_{i=0}^{\ell} \binom{n}{i} = \binom{1}{0} +\binom{1}{1} = 2$. For $n>1$, we assume for the induction hypothesis that all union-closed families $\mathcal{A}'$ with $|U(\mathcal{A}')| < n$ have $|\mathcal{A}'| \leq \sum_{i=0}^{\ell(\mathcal{A}')} \binom{|U(\mathcal{A}')|}{i}$, with equality if and only if $\mathcal{A}'=\bigcup_{i=0}^{\ell(\mathcal{A}')}\binom{U(\mathcal{A}')}{|U(\mathcal{A}')|-i}$. We let $\mathcal{A}$ be any union-closed family of sets with universe size $n>1$ and length $\ell < n$. (If $\ell=n$, then the theorem is satisfied, as $|\mathcal{A}| \leq \sum_{i=0}^{\ell} \binom{n}{i} = 2^n$ with equality if and only if $\mathcal{A}= \bigcup_{i=0}^{\ell}\binom{[n]}{n-i} = 2^{[n]}$.) For any $x \in [n]$, we define the following families:

\vspace{0.075cm}

\[\ax = \{A \in \mathcal{A} \ | \ x \in A \}\textrm{;}\]
\[\hax = \{A \setminus \{x\} \ | \ A \in \ax\}\textrm{;}\]
\[\abx=\{A \in \mathcal{A} \ | \ x \not \in A \}\textrm{.}\]

\vspace{0.275cm}

\noindent We observe that $\mathcal{A} = \ax \cup \abx$, $\ax \cap \abx = \emptyset$, and $|\ax| = |\hax|$ together imply that $|\mathcal{A}| = |\hax| + |\abx|$. Because no member set of $\hax$ contains $x$, and $[n] \in \ax$ implies that $[n] \setminus \{x\} \in \hax$, we have that $|\uhax| = n-1$. Similarly, because no member set of $\abx$ contains $x$, we have that $|\abx| \leq n-1$. Now, since $\ax$ is a subfamily of $\mathcal{A}$ and $\lhax = \lax$, we have that $\lhax \leq \ell$. Let $\ell(\{\emptyset\}) = 0$. Noting that $\abx$ is a subfamily of $\mathcal{A}$, that $[n]$ belongs to any chain of maximum size in $\mathcal{A}$, and that $[n]$ does not belong to $\abx$, we also have that either $\ell=0$ or $\lanx \leq \ell-1$. Since $\hax$ is itself union-closed, we have by the induction hypothesis that

\medskip

\[|\hax| \leq \sum_{i=0}^{\lhax} \binom{|\uhax|}{i} = \sum_{i=0}^{\lhax} \binom{n-1}{i} \leq \sum_{i=0}^{\ell}\binom{n-1}{i}\textrm{,}\]

\bigskip

\noindent If $\abx = \emptyset$, then $|\mathcal{A}|=|\hax| \leq \sum_{i=0}^{\ell}\binom{n-1}{i} \leq \sum_{i=0}^{\ell} \binom{n}{i}$. If $\abx = \{\emptyset\}$, then  $\ell > 0$ and $|\mathcal{A}|=|\hax|+1 \leq 1 + \sum_{i=0}^{\ell}\binom{n-1}{i} \leq \sum_{i=0}^{\ell} \binom{n}{i}$. Else, $\abx \neq \emptyset$ and $\abx \neq \{\emptyset\}$, and $\abx$ is union-closed as well as $\hax$. In this case, $\ell$ is again positive, and we apply the induction hypothesis to $\abx$ in order to obtain that

\medskip
\smallskip

\[|\abx| \leq \sum_{i=0}^{\lanx} \binom{|\uabx|}{i}\leq \sum_{i=0}^{\lanx} \binom{n-1}{i} \leq \sum_{i=0}^{\ell-1}\binom{n-1}{i}\textrm{.}\]

\bigskip
\smallskip

\noindent Therefore, if $\abx \neq \emptyset$ and $\abx \neq \{\emptyset\}$, then we have that

\vspace{-0.1cm}

\begin{align*}
&\hspace{-0.4375cm}|\mathcal{A}| = |\hax| + |\abx| \\ 
& \leq \sum_{i=0}^{\ell} \binom{n-1}{i} + \sum_{i=0}^{\ell-1} \binom{n-1}{i} \\
& = \binom{n-1}{0} + \sum_{i=1}^{\ell} \Bigr(\binom{n-1}{i} + \binom{n-1}{i-1}\Bigr) \\
& = 1 + \sum_{i=1}^{\ell} \binom{n}{i} \\
& = \sum_{i=0}^{\ell} \binom{n}{i}\textrm{.}
\end{align*} 

\bigskip

\noindent (Here, the equation $\sum_{i=0}^{\ell} \binom{n}{i} = \sum_{i=0}^{\ell} \binom{n-1}{i} + \sum_{i=0}^{\ell-1} \binom{n-1}{i}$ corresponds to case $m=1$ of the identity $\sum_{i=0}^k \binom{n}{i} = \sum_{i=0}^m \binom{m}{i} \sum_{j=0}^{k-i} \binom{n-m}{j}$ for all integers  $0 \leq k \leq n$, where $0 \leq m \leq n$.) This completes the proof that $|\mathcal{A}| \leq \sum_{i=0}^{\ell} \binom{n}{i}$. Next, we show that  $|\mathcal{A}| = \sum_{i=0}^{\ell} \binom{n}{i}$ if and only if $\mathcal{A} = \bigcup_{i=0}^{\ell}\binom{[n]}{n-i}$. That $|\mathcal{A}| = \sum_{i=0}^{\ell} \binom{n}{i}$ when $\mathcal{A} = \bigcup_{i=0}^{\ell}\binom{[n]}{n-i}$ is clear from symmetry of the binomial coefficient. It remains to show that the bound is only sharp if $\mathcal{A} = \bigcup_{i=0}^{\ell}\binom{[n]}{n-i}$, for which we now assume that $|\mathcal{A}| = \sum_{i=0}^{\ell} \binom{n}{i}$. If $\ell=0$, then $|\mathcal{A}| = \sum_{i=0}^{\ell}\binom{n}{i} = \binom{n}{0}=1$, implying that $\mathcal{A} = \{[n]\} = \binom{[n]}{n} = \bigcup_{i=0}^{\ell} \binom{[n]}{n-i}$. Else, $\ell > 0$, and we have that $\sum_{i=0}^{\ell} \binom{n}{i} = \sum_{i=0}^{\ell} \binom{n-1}{i} + \sum_{i=0}^{\ell-1} \binom{n-1}{i}$. Noting that $|\mathcal{A}| = |\hax|+|\abx|$ and $|\hax| \leq \sum_{i=0}^{\ell} \binom{n-1}{i}$, we then have that $|\abx| \geq \sum_{i=0}^{\ell-1} \binom{n-1}{i}$. It again follows that $\abx$ is union-closed. (Otherwise, $\abx = \emptyset$ and $|\abx|=0 \geq \sum_{i=0}^{\ell-1}\binom{n-1}{i} \geq 1$, a contradiction, or $\abx = \{\emptyset\}$ and $|\abx| = 1 \geq \sum_{i=0}^{\ell-1}\binom{n-1}{i}$ implies that $\ell=1$, making $\mathcal{A} = \{[n], \emptyset\}$, which then implies that $n=1$, again a contradiction). Therefore, $|\abx| \leq \sum_{i=0}^{\ell-1} \binom{n-1}{i}$, which implies that both $|\hax| = \sum_{i=0}^{\ell} \binom{n-1}{i}$ and $|\abx| = \sum_{i=0}^{\ell-1} \binom{n-1}{i}$. It follows that $\uhax = [n] \setminus \{x\}$ with $\lhax = \ell$. (If not, then $|\hax| = \sum_{i=0}^{\lhax} \binom{|\uhax|}{i} \leq \sum_{i=0}^{\ell} \binom{n-2}{i}$ or $|\hax| = \sum_{i=0}^{\lhax} \binom{|\uhax|}{i} \leq \sum_{i=0}^{\ell-1} \binom{n-1}{i}$, a contradiction.) Therefore, $\hax = \bigcup_{i=0}^{\ell} \binom{[n] \setminus \{x\}}{(n-1)-i}$ by the induction hypothesis. Similarly, it also follows that $\uabx = [n] \setminus \{x\}$ with $\lanx = \ell-1$. (When $\ell=1$, having $x \not \in A$ for any $A \in \abx$ implies that $|\abx|=1$, making $\ell(\abx)=0=\ell-1$.) Applying the induction hypothesis now to $\abx$, we obtain that $\abx = \bigcup_{i=0}^{\ell-1} \binom{[n] \setminus \{x\}}{(n-1)-i}$. We thus have that

\vspace{-0.1cm}

\begin{align*}
&\hspace{-0.3cm}\mathcal{A} = \ax \cup \abx = \{A \cup \{x\} \ | \ A \in \hax\} \cup \abx \\
& = \bigcup_{i=0}^{\ell} \Bigr\{A \cup \{x\} \ \Bigr | \ A \in \binom{[n] \setminus \{x\}}{n-i-1} \Bigr\} \cup \bigcup_{i=0}^{\ell-1} \binom{[n]\setminus\{x\}}{n-i-1} \\
& = \{[n]\} \cup \bigcup_{i=1}^{\ell} \Bigr(\Bigr\{A \cup \{x\} \ \Bigr | \ A \in \binom{[n] \setminus \{x\}}{n-i-1} \Bigr\} \cup \binom{[n]\setminus \{x\}}{n-i}\Bigr) \\
& = \binom{[n]}{n} \cup \bigcup_{i=1}^{\ell} \binom{[n]}{n-i} \\
& = \bigcup_{i=0}^{\ell} \binom{[n]}{n-i}\textrm{.}
\end{align*} 

\medskip

\noindent Therefore, $|\mathcal{A}| = \sum_{i=1}^{\ell}\binom{n}{i}$ implies that $\mathcal{A} = \bigcup_{i=0}^{\ell}\binom{[n]}{n-i}$, completing the proof of Theorem 1.

\bigskip
\smallskip

\noindent \textbf{Corollary 2.1.} \textit{For any union-closed family} $\mathcal{A}$\textit{, there is an element from its universe} $[n]$ \textit{that is in at most} $\sum_{i=0}^{\ell} \binom{n-1}{i}$\textit{ of its member sets.}

\medskip

\noindent \textit{Proof.} Because $|\mathcal{A}| \leq \sum_{i=0}^{\ell}\binom{n}{i}$, it holds that $|\mathcal{A} \setminus \bigcup_{i=0}^{\ell} \binom{[n]}{n-i}| \leq |\bigcup_{i=0}^{\ell} \binom{[n]}{n-i} \setminus \mathcal{A}|$. We also have that $|X| < |Y|$ for any $X \in \mathcal{A} \setminus \bigcup_{i=0}^{\ell} \binom{[n]}{n-i}$ and $Y \in \bigcup_{i=0}^{\ell} \binom{[n]}{n-i} \setminus \mathcal{A}$. Together these observations imply that $\sum_{A \in \mathcal{A}}|A| \leq \sum_{A \in \bigcup_{i=0}^{\ell} \binom{[n]}{n-i}}|A| = \sum_{i=0}^{\ell} (n-i)\binom{n}{i}$. Because $\sum_{A \in \mathcal{A}}|A| = \sum_{x \in [n]}|\ax|$, we also have that $\sum_{x \in [n]}|\ax| \leq \sum_{i=0}^{\ell}(n-i)\binom{n}{i}$. It follows that $|\mathcal{A}_y| \leq \sum_{i=0}^{\ell}(1-\frac{i}{n})\binom{n}{i}$ for some $y \in [n]$. (If not, then $|\ax| > \sum_{i=0}^{\ell}(1-\frac{i}{n})\binom{n}{i}$ for all $x \in [n]$, implying that $\sum_{x \in [n]}|\ax| > n(\sum_{i=0}^{\ell}(1-\frac{i}{n})\binom{n}{i})=\sum_{i=0}^{\ell}(n-i)\binom{n}{i}$, a contradiction.) Finally, considering that $\sum_{i=0}^{\ell}(1-\frac{i}{n})\binom{n}{i}$ is equal to the size of $\ax$ for $\mathcal{A}=\bigcup_{i=0}^{\ell} \binom{[n]}{n-i}$ and any $x \in [n]$, we have by double counting that $\sum_{i=0}^{\ell}(1-\frac{i}{n})\binom{n}{i} = \sum_{i=0}^{\ell} \binom{n-1}{n-i-1}$, making $|\mathcal{A}_y| \leq \sum_{i=0}^{\ell} \binom{n-1}{i}$.

\vspace{0.325cm}

\section*{3. Proof of Theorem 2}

\noindent Let $\mathcal{A}$ be any union-closed family with universe $[n]$ and length $\ell$, and consider any chain $\mathcal{C}=\{C_1, \cdots, C_{\ell+1}\}$ in $\mathcal{A}$ of maximum size, where $1 \leq i < j \leq \ell+1$ implies that $C_j \subsetneq C_i$ without loss of generality. Because every member set of $\mathcal{A}$ is a subset of $[n]$, and $[n]$ must belong to $\mathcal{A}$, we have that $C_1=[n]$. Let $[0] = \emptyset$. For each $i \in [\ell]$, we define

\[\mathcal{D}_i = \ \Bigr\{X \setminus (C_i \setminus C_{i+1}) \ \Bigr | \ X \in \mathcal{C}_i\Bigr\}\textrm{,}\]

\vspace{0.3875cm}

\noindent where

\vspace{-0.3875cm}

\noindent \[\mathcal{C}_i = \Bigr\{X \in \mathcal{A} \ \Bigr| \ (C_i \setminus C_{i+1} \subseteq X) \land (X \cap \bigcup_{j \in [i-1]} C_j \setminus C_{j+1} = \emptyset)\Bigr\}\textrm{.}\]

\bigskip

\noindent \textbf{Theorem 3.1.} \textit{For any} $i \in [\ell]$\textit{, if} $C_{i+1} \neq \emptyset$\textit{, then} $\mathcal{D}_i$ \textit{is union-closed.}

\medskip

\noindent \textit{Proof.} Let $i$ be any element of $[\ell]$ such that $C_{i+1} \neq \emptyset$. Assuming any $X_1,X_2 \in \mathcal{C}_i$, we first have that $(C_i \setminus C_{i+1} \subseteq X_1) \ \land \ (C_i \setminus C_{i+1} \subseteq X_2) \implies C_i \setminus C_{i+1} \subseteq X_1 \cup X_2$ and $(X_1 \cap \bigcup_{j \in [i-1]} C_j \setminus C_{j+1} = \emptyset) \ \land \ (X_2 \cap \bigcup_{j \in [i-1]} C_j \setminus C_{j+1} = \emptyset) \implies (X_1 \cup X_2) \cap \bigcup_{j \in [i-1]} C_j \setminus C_{j+1} = \emptyset$. Hence, $X_1 \cup X_2$ satisfies both conditions for membership in $\mathcal{C}_i$. Because $C_i \in \mathcal{C}_i$ and $C_i \neq \emptyset$, it then follows that $\mathcal{C}_i$ is union-closed. Now, assuming any $Y_1,Y_2 \in \mathcal{D}_i$, we have that $Y_1 \cup (C_i \setminus C_{i+1}) \in \mathcal{C}_i$ and $Y_2 \cup (C_i \setminus C_{i+1}) \in \mathcal{C}_i$. Since $\mathcal{C}_i$ is union-closed, we then also have that $(Y_1 \cup (C_i \setminus C_{i+1})) \cup (Y_2 \cup (C_i \setminus C_{i+1}))= (Y_1 \cup Y_2) \cup (C_i \setminus C_{i+1}) \in \mathcal{C}_i$. It follows that $((Y_1 \cup Y_2) \cup (C_i \setminus C_{i+1})) \setminus (C_i \setminus C_{i+1}) = Y_1 \cup Y_2 \in \mathcal{D}_i$. Because $C_{i+1} \in \mathcal{D}_i$ and $C_{i+1} \neq \emptyset$, it then follows that $\mathcal{D}_i$ is also union-closed, proving Theorem 3.1.

\bigskip
\medskip

\noindent \textbf{Theorem 3.2.} \vspace{-0.55cm} \[\hspace{-0.5cm}|\mathcal{A}| = 1+ \sum_{i=1}^{\ell} |\mathcal{D}_i|\textit{.}\]

\medskip

\noindent \textit{Proof.} We first have that $\mathcal{A} \setminus \{C_{\ell+1}\} = \bigcup_{i \in [\ell]} \normalfont \mathcal{C}_i$. (If not, then there exists some $X \in \mathcal{A} \setminus \{C_{\ell+1}\}$ such that either $X \subsetneq C_{\ell+1}$ or $\exists i \in [\ell] \ | \ (\emptyset \subsetneq (C_i \setminus C_{i+1}) \cap X \subsetneq C_i \setminus C_{i+1} \ \land \ \forall j \in [i-1] \ | \ (C_j \setminus C_{j+1}) \cap X = \emptyset)$, and it follows that either $\{C_1, \cdots, C_{\ell + 1}, X\}$ or $\{C_1, \cdots, C_i, C_{i+1} \cup X, C_{i+1}, \cdots, C_{\ell + 1}\}$ is a chain in $\mathcal{A}$ of size $\ell+2$, contradicting the definition of $\ell$.) Next, we observe that $\mathcal{C}_i \cap \mathcal{C}_j = \emptyset$ for any distinct $i$ and $j$ in $[\ell]$. (Otherwise, there exist distinct $i$ and $j$ in $[\ell]$ with at least one member set $X$ in $\mathcal{C}_i \cap \mathcal{C}_j$, where $j < i$ without loss of generality, and $X \in \mathcal{C}_j$ implies that $C_j \setminus C_{j+1} \subseteq X$, while $X \in \mathcal{C}_i$ and $j \in [i-1]$ together imply that $X \cap (C_j \setminus C_{j+1}) = \emptyset$, a contradiction.) It follows that $|\mathcal{A} \setminus \{C_{\ell+1}\}| = \sum_{i=1}^{\ell}|\mathcal{C}_i|$, making $|\mathcal{A}| = 1 + \sum_{i=1}^{\ell}|\mathcal{C}_i|$. Finally, because $C_i \setminus C_{i+1} \subseteq X$ for any $i \in [\ell]$ and any $X \in \mathcal{C}_i$, we have that $|\mathcal{D}_i|=|\mathcal{C}_i|$ for each $i$. It then follows that $|\mathcal{A}| = 1+\sum_{i=1}^{\ell}|\mathcal{D}_i|$, proving Theorem 3.2.

\bigskip

\noindent \textbf{Lemma 3.3.} $\ \ |U(\mathcal{D}_i)| \leq n-i$ \textit{and} $\ell(\mathcal{D}_i) \leq \ell$ \textit{for all} $i \in [\ell]$.

\medskip

\noindent \textit{Proof.} For all $i \in [\ell]$, $|U(\mathcal{D}_i)| \leq n-i$ follows from having $X \cap \bigcup_{j \in [i]} C_j \setminus C_{j+1} = \emptyset$ for any $X \in \mathcal{D}_i$, and $\ell(\mathcal{D}_i) \leq \ell$ follows from $\ell(\mathcal{C}_i) = \ell(\mathcal{D}_i)$ and $\mathcal{C}_i$ being a subfamily of $\mathcal{A}$.

\bigskip

\noindent We now define $\Theta \colon {\mathbb{Z}^3_{\geq 0}} \to \mathbb{Z}_{\geq 0}$ such that \[\Theta(x,y,z) = \frac{x^z-1}{x-1} + 2^y(1-2^{-x})^z\textrm{.}\]

\medskip

\noindent We prove by induction that, for any nonnegative integer $p$, $|\mathcal{A}| \leq \Theta(\ell, n , p)$ for all union-closed families $\mathcal{A}$ with universe $[n]$ and length $\ell$. For the base case $p=0$, we have that $|\mathcal{A}| \leq \Theta(\ell, n , p) = 2^n$ for any such $\mathcal{A}$. For $p > 0$, we again let $\mathcal{A}$ be any union-closed family with universe $[n]$ and length $\ell$, and we assume for the induction hypothesis that any union-closed family $\mathcal{A}'$ with universe of size $n'$ and length $\ell'$ satisfies $|\mathcal{A}'| \leq \Theta(\ell',n', p-1)$. By Theorem 3.1, we apply this hypothesis to all families $\mathcal{D}_i \colon i \in [\ell]$ that have $C_{i+1} \neq \emptyset$ in order to obtain that $|\mathcal{D}_i| \leq \Theta(\ell(\mathcal{D}_i),|U(\mathcal{D}_i)|,p-1)$ for all such $i$. (If $C_{i+1} = \emptyset$, then we may directly compute that $|\mathcal{D}_i| = \Theta(\ell(\mathcal{D}_i),|U(\mathcal{D}_i)|,p-1)=\Theta(0,0,p-1)=1$.) Noting that $\Theta(x,y,z)$ is always increasing with respect to each of the variables $x$ and $y$, it then follows from Lemma 3.3 that $|\mathcal{D}_i| \leq \Theta(\ell,n-i,p-1)$ for all $i \in [\ell]$. Recall that by Theorem 3.2, $|\mathcal{A}| = 1+\sum_{i=1}^{\ell} |\mathcal{D}_i|$. We  substitute $\Theta(\ell, n-i, p-1)$ into this equation for every $|\mathcal{D}_i|$ to obtain an upper bound for the size of $\mathcal{A}$ as follows:

\vspace{-0.14cm}

\begin{align*}
&\hspace{0.715cm} |\mathcal{A}| = 1+\sum_{i=1}^{\ell}|\mathcal{D}_i|
\\
&\hspace{0.575cm}\hspace{0.575cm} \leq 1+\sum_{i=1}^{\ell}\Theta(\ell,n-i,p-1)
\\
&\hspace{0.575cm}\hspace{0.575cm} = 1+\sum_{i=1}^{\ell}\Bigr( \frac{\ell^{p-1}-1}{\ell-1} + 2^{n-i}(1-2^{-\ell})^{p-1}\Bigr)
\\ 
&\hspace{0.575cm}\hspace{0.575cm}= 1+ \ell\Bigr(\frac{\ell^{p-1}-1}{\ell-1}\Bigr)+(1-2^{-\ell})^{p-1}\sum_{i=1}^{\ell}2^{n-i}
\\
&\hspace{0.575cm}\hspace{0.575cm}=\frac{\ell^p-1}{\ell-1}+(1-2^{-\ell})^{p-1}(2^n-2^{n-\ell})
\\
&\hspace{0.575cm}\hspace{0.575cm}=\frac{\ell^p-1}{\ell-1}+2^n(1-2^{-\ell})^p
\\
&\hspace{0.575cm}\hspace{0.575cm}=\Theta(\ell, n , p)\textrm{.}
\end{align*}

\vspace{0.325cm}

\noindent This resolves the induction step. It follows that $|\mathcal{A}| \leq  \min_{p \in \mathbb{Z}_{\geq 0}}\{\Theta(\ell, n , p)\}$ for any union-closed family $\mathcal{A}$ with universe $[n]$ and length $\ell$.

\medskip

\noindent For all integers $1 \leq k \leq n$, we set $\mathcal{A}$ equal to $\bigcup_{i=0}^{k} \binom{[n]}{n-i}$ to obtain that $\sum_{i=0}^k \binom{n}{i} \leq \min_{p \in \mathbb{Z}_{\geq 0}}\{\Theta(k, n , p)\}$. As $\hat{p} \in \mathbb{Z}_{\geq 0}$, this completes the proof of Theorem 2.

\bigskip

To demonstrate that $\hat{p}$ is optimal, we compute that $\Theta(k,n,p+1) \leq \Theta(k,n,p)$ if and only if \[p \leq \log_{\frac{k}{1-2^{-k}}}(2^{n-k}) = \frac{n-k}{\log_2(k/(1-2^{-k}))}\textrm{.}\]

\medskip

Thus, $\Theta$ is decreasing with respect to $p$ if and only if $p \leq \lfloor \frac{n-k}{\log_2(k/(1-2^{-k}))} \rfloor$, and it follows that \[\Theta \Bigr(k,n,\Bigr\lfloor \frac{n-k}{\log_2(k/(1-2^{-k}))} \Bigr\rfloor+1\Bigr) = \Theta(k,n,\hat{p}) = \min_{p \in \mathbb{Z}_{\geq 0}} \{ \Theta(k,n,p) \}\textrm{.}\]

\end{document}